\documentclass[12pt,a4paper]{article}
\usepackage{amsmath,amssymb,amsfonts}
\def\Bbb{\mathbb}

\title{\bf On Dedekind sums with equal values}

\author{Kurt Girstmair\\ \small
Institut f\"ur Mathematik, Universit\"at Innsbruck   \\
\small Technikerstr. 13/7, A-6020 Innsbruck, Austria \\
\small Kurt.Girstmair@uibk.ac.at}

\date{}

\normalsize

\makeatletter
\let\@@maketitle=\maketitle
\def\maketitle{\def\thispagestyle##1{\relax}\@@maketitle}
\makeatother
\newtheorem{theorem}{Theorem}

\newtheorem{corollary}{Corollary}

\def\BE{\begin{equation}}
\def\EE{\end{equation}}
\def\BD{\begin{displaymath}}
\def\ED{\end{displaymath}}
\def\BA{\begin{array}}
\def\EA{\end{array}}
\def\BEA{\begin{eqnarray*}}
\def\EEA{\end{eqnarray*}}
\def\BI{\bibitem}

\def\Z{\Bbb Z}

\def\R{\Bbb R}

\def\phi{\varphi}
\def\EPS{\varepsilon}

\def\CMOD#1#2#3{#1 \equiv #2 \: \mbox{mod}\: #3}

\def\MB{\mbox}
\def\LD{\ldots}

\def\DIV{\,|\,}
\def\NDIV{\, \nmid \,}
\def\LS#1#2{\left(\frac{#1}{#2}\right)}

\def\BQ{``}
\def\EQ{'' }

\def\STOP{\hfill$\Box$}

\def\DED{Dedekind }

\begin{document}
\maketitle

\begin{abstract}

\noindent
\DED sums $s(m,n)$ occur in many fields of mathematics. Since $s(m_1,n)=s(m_2,n)$ if $\CMOD{m_1}{m_2}n$,
it is natural to ask which of the \DED sums $s(m,n)$, $0\le m<n$, take equal values. So far no
simple criterion is known by which the equality of $s(m_1,n)$ and $s(m_2,n)$ could be decided. In this note
we show how to obtain non-obvious examples of equal \DED sums. We consider two cases which mark the extreme
possibilities for the argument $n$, namely, $n$ a prime power and $n$ square-free.
Whereas we can give a partial overview of equal \DED sums in the prime power case, such an overview seems
to be much more difficult to obtain in the square-free case.

\end{abstract}

\section*{Introduction}

Let $n$ be a positive integer and $m\in \Z$, $(m,n)=1$. The classical {\em \DED sum} $s(m,n)$ is defined by
\BE
\label{0.2}
   s(m,n)=\sum_{k=1}^n ((k/n))((mk/n))
\EE
where $((\LD))$ is the \BQ sawtooth function\EQ defined by
\BD
  ((t))=\left\{\begin{array}{ll}
                 t-\lfloor t\rfloor-1/2 & \MB{ if } t\in\R\smallsetminus \Z; \\
                 0 & \MB{ if } t\in \Z
               \end{array}\right.
\ED
(see, for instance, \cite[p. 1]{RaGr}).

\DED sums have quite a number of interesting applications in
analytic number theory (modular forms), algebraic number theory (class numbers),
lattice point problems, topology and algebraic geometry (see, for instance, \cite{{Ap}, {At}, {BeRo}, {Me}, {RaGr}, {Ur}}).
Moreover, various properties of these sums have been studied by several authors
(see, for instance, \cite{{Ba}, {Br}, {Gi}, {Hi}, {Kn}, {Va}, {Zh}}).

In the present setting it is more convenient to work with
\BD
 S(m,n)=12s(m,n)
\ED instead. Observe that $S(m_1,n)=S(m_2,n)$ if $\CMOD{m_1}{m_2}n$, so one often considers
only arguments $m$ in the range $0\le m< n$.

If we fix $n$, we may ask which of the \DED sums $S(m,n)$, $0\le m< n$, $(m,n)=1$, take equal values. In
the paper \cite{Ja} it was shown that $S(m_1,n)=S(m_2,n)$ only if
\BE
\label{0.4}
      (m_1-m_2)(m_1m_2-1)\equiv 0 \MB{ mod } n.
\EE
This condition, however, is not sufficient for the equality of $S(m_1,n)$ and $S(m_2,n)$. Indeed, the
condition is necessary and sufficient for
\BD
  S(m_1,n)-S(m_2,n)\in\Z
\ED
(see \cite{Gi2}). It seems that a simple necessary and sufficient condition for the equality of $S(m_1,n)$ and $S(m_2,n)$
is currently out of reach.

Suppose, for the moment, that $n=p_1p_2\LD p_r$ is square-free (so $p_1,\LD,p_r$ are distinct primes) and that
$m_1$ is given. It is known that
the number of integers $m_2$, $0\le m_2<n$, $(m_2,n)=1$, such that $S(m_1,n)-S(m_2,n)\in\Z$ is $\le 2^r$
(see \cite[Th. 3]{Gi2}). In particular, there are at most $2^r$ numbers $m_2$ in this range with $S(m_1,n)=S(m_2,n)$.
Accordingly, the \DED sums $S(m,n)$, $0\le m<n$, $(m,n)=1$, take at least
\BD
\prod_{j=1}^r\frac{p_j-1}2
\ED
distinct values. So there are, as a rule, plenty of values $S(m,n)$ that must be distinguished.

In view of this situation, it may be worthwhile exhibiting {\em series} of equal \DED sums.
To this end we apply two
theorems from the literature (one of Rademacher and one of our own). Whereas the first theorem
gives insight into the case of powers $n=l^k$, $k\ge 2$ (so it comprises, in particular, the case of prime powers),
the second one supplies examples of equal \DED sums for square-free numbers $n$. In the prime power case $n=p^k$ we
obtain a partial overview of the equalities $S(m_1,n)=S(m_2,n)$ that can occur in this situation.
In the square-free case such an overview seems to be much more difficult to obtain.

In all of these examples we distinguish between {\em obvious} equality and {\em non-obvious} equality. Indeed,
it is almost obvious from (\ref{0.2}) that $S(m_,n)=S(m^*,n)$, where $m^*$ is a {\em multiplicative inverse} of $m$
mod $n$, i.e., $\CMOD{mm^*}1n$ (see \cite[p. 26]{RaGr}). This case of equality is addressed as the obvious case,
whereas all other cases are considered as non-obvious.

\section*{1. The power case}

In the paper \cite{Ra}, Rademacher enunciated his Satz 15 in a way which does not immediately show its applicability
to equal \DED sums. Here we note a slightly weaker version of Rademacher's result, which, however, obviously
produces examples of equal \DED sums, namely,

\begin{theorem} 
\label{t1}
Let $d$ and $n$ be positive integers and $m\in \Z$, $(m,n)=1$. Let $\EPS\in\{\pm 1\}$.
Then
\BE
\label{1.4}
  S(\EPS+dnm,dn^2)=\EPS\left(\frac 2{dn^2}+d-3\right).
\EE

\end{theorem} 

\medskip
\noindent
Rademacher proved his theorem by means of invariants of binary quadratic forms. In Section 3 we shall give a
proof of Theorem \ref{t1} by means of the three-term relation for \DED sums (which is due to Rademacher and Dieter).
Therefore, this proof is, in some sense, more at home in the setting of \DED sums
than Rademacher's proof. Moreover, our proof may serve as a model for the proof
of Theorem \ref{t2}, which is the basis of Section 2.

In the setting of Theorem \ref{t1}, let $m$ run through all integers $0\le m<n$, $(m,n)=1$.
The theorem says that each of the \DED sums $S(\EPS+dnm,dn^2)$ takes the same value.
So whenever $\phi(n)>2$, there must be non-obvious cases of equal values among them.

\medskip
\noindent
{\em Example.} Let $d=8$, $n=5$. Then $S(1+40m,200)=1/100+5=5.01$ for $m=1,2,3,4$. Here the equality $S(41,200)=S(81,200)$
is non-obvious in  the above sense, whereas $S(41,100)=S(161,100)$ is obvious.

\medskip
In the case of powers $n=l^k$, Rademacher's theorem and Theorem \ref{t1} give the same, namely,

\begin{corollary} 
\label{c1}
Let $l$, $k$, $r$, $q$ be positive integers and $r\le k$.
Suppose that $q\DIV l^{k-r}$ and $l\NDIV q$. Suppose, further, that one of the following holds:

\rm{(a)} $r\ge k/2$;

\rm{(b)} $r<k/2$ and $l^{k-2r}\DIV q^2$.

\noindent Let $\EPS\in\{\pm 1\}$.
If $m$ is an integer such that $(m,l^{k-r}/q)=1$, then
\BE
\label{1.6}
  S(\EPS+l^rqm,l^k)=\EPS\left(\frac 2{l^k}+l^{2r-k}q^2-3\right).
\EE

\end{corollary} 

\noindent
{\em Proof of Corollary \rm{\ref{c1}}.}
In the setting of this corollary, put
$n=l^{k-r}/q$ (a positive integer) and $d=l^ {2r-k}q^2$. In the case of assumption (a), $2r-k\ge 0$, so
$d$ is a positive integer. This is also true in the case of assumption (b).
Then $dn=l^rq$ and $dn^2=l^k$. Accordingly, Theorem \ref{t1} gives (\ref{1.6})
for each integer $m$ with $(m, l^{k-r}/q)=1$.\STOP

\medskip
\noindent
{\em Examples.} 1. Let $l=6$, $k=4$, $r=2$ and $q=4$ (so assumption (a) holds).
Then $l^{k-r}/q=36/4=9$ and $l^rq=144$. For $m=1,2,4,5,7,8$
we have $S(1+144m, 1296)=2/1296+16-3=2/1296+13$. Again, the equality of the values is non-obvious for $m=1,2,4$.

2. The case (b) is illustrated by the following example. Let $l=12$, $k=3$, $r=1$ and $q=6$.
Here $12\, (=l^{k-2r})$ divides $36\,(=q^2)$. We obtain
$S(1+72m, 1728)=2/1728+36/12-3=2/1728$ for $m$ relatively prime to $l^{k-r}/q=24$. The equality of the values is
non-obvious for $m=1,5,7,11$.

\medskip
\noindent
If $l=p$ is a prime number, only the case $q=1$ is possible since $p\NDIV q$ and $q\DIV p^{k-r}$.
Accordingly, only case (a) of Corollary \ref{c1} applies here. In this case, however, we have much more information
about the equality of \DED sums. Indeed, suppose that $(m_1,p)=1$ and $m_1\not\equiv \pm 1$ mod $p$. If $m_2$
satisfies (\ref{0.4}), we either have $\CMOD{m_1}{m_2}p$ or $\CMOD{m_1m_2}1p$. Each of these cases excludes the
other because $\CMOD{m_1^2}1p$ implies $\CMOD{m_1}{\pm 1}p$. So we are left
with $\CMOD{m_1}{m_2}{p^k}$ or $\CMOD{m_1m_2}1{p^k}$. In other words, the assumption
$m_1\not\equiv\pm 1$ mod $p$ allows only obvious equalities $S(m_1,p^k)=S(m_2,p^k)$.

Hence we have to consider only numbers $m_1$ of the form $m_1=\EPS+p^rm$, $\EPS\in\{\pm 1\}$, $p\NDIV m$, $1\le r\le k$.
Here we obtain a complete overview of equal values of \DED sums if we assume
$r\ge k/2$. Let this assumption hold. In \cite[equ. (9)]{Gi3} we have shown that $S(m_1,p^k)-S(m_2,p^k)\in\Z$ only if
$m_2=\EPS+p^{r'}m'$, $r'\ge k/2$, $p\NDIV m'$. By Corollary \ref{c1}, we have
\BD
  S(m_1,p^k)=\EPS(2/p^k+p^{2r-k}-3) \MB{ and } S(m_2,p^k)=\EPS(2/p^k+p^{2r'-k}-3).
\ED
So these values are equal if, and only if, $r=r'$. Altogether, we obtain

\begin{corollary} 
\label{c2}
Let $p$ be a prime number and $k$, $r$ positive integers with $k/2\le r\le k$.
Let $\EPS\in\{\pm 1\}$ and $m_1=\EPS+p^rm$ with an integer $m$, $p\NDIV m$.
For $m_2\in\Z$, $p\NDIV m_2$, we have $S(m_2,p^k)=S(m_1,p^k)$ if, and only if,
$m_2=\EPS+p^rm'$, $m'\in\Z$, $p\NDIV m'$. In this case
\BD
\label{1.8}
  S(m_1,p^k)=S(m_2,p^k)=\EPS\left(\frac 2{p^k}+p^{2r-k}-3\right).
\ED

\end{corollary} 

\noindent
{\em Remark.} In the case $1\le r<k/2$ there is apparently no result like Corollary \ref{c2}. Consider, for instance,
$p=3$, $k=5$, $r=2$. There are 18 values $1+9m$, $3\NDIV m$, in the range $0\le 1+9m<243$. It suffices to consider
$m=1,2,4,10,11,13,14,16, 17$, since the remaining values $1+9m$ arise from these as multiplicative inverses mod $243$.
The corresponding \DED sums have the form
\BD
   S(1+9m,243)=\frac{83}{243}+z,
\ED
with $z\in\{-27,-19,-11,-3,\,5,\,13,\, 21\}$ (so all of these \DED sums have the same fractional part).
The only equal values among these occur for $m=4$ and $m'=14$ (with $z=5$),
and for $m=11$ and $m'=16$ (with $z=-11$).

\medskip
\noindent
If we apply the above considerations (in particular, Corollary \ref{c2}) to the case $k=2$, $r=1$,
we obtain

\begin{corollary} 
\label{c3}
Let $p$ be a prime number and $\EPS\in\{\pm 1\}$.
Then all values $S(\EPS+pm,p^2$), $m=1,\LD,p-1$, are equal, namely,
\BD
\label{1.10}
  S(\EPS+pm,p^2)=\EPS\left(\frac 2{p^2}-2\right).
\ED
If $p\ge 5$, we have, thus, non-obvious equalities
$S(m_1,p^2)=S(m_2,p^2)$ for $m_1, m_2\in\{1+pm;m=1,\LD,(p-1)/2\}$, $m_1\ne m_2$.
All other non-obvious equalities
$S(m_1,p^2)=S(m_2,p^2)$, $m_1,m_2\in\{1,\LD,p^2-1\}$, $p\NDIV m_1,m_2$,
arise from these by
transition to multiplicative inverses mod $p^2$.

\end{corollary} 

\section*{2. The square-free case}

Many examples of non-obvious equalities in the square-free case arise from

\begin{theorem} 
\label{t2}
Let $n$ be a positive integer and $m\in \Z$, $(m,n)=1$.
As above, let $m^*\in\Z$ denote an inverse of $m$ mod $n$, i.e., $\CMOD{mm^*}1n$.
Let $t$ be a positive integer with $\CMOD t{m-m^*}n$.
Then
\BE
\label{2.2}
  S(1+mt,nt)=\frac 2{nt}+\frac tn-3.
\EE

\end{theorem} 

\noindent
For a proof of Theorem \ref{t2}, see \cite{Gi4}. In Section 3 we briefly show how to adapt the proof of
Theorem \ref{t1} in order to obtain a proof of Theorem \ref{t2}.  In what follows,
$\LS qp$ denotes the Legendre symbol for an integer $q$ and a prime $p$.

\begin{corollary} 
\label{c4}
Let $t$ be a positive integer and $t^2+4=qk^2$, where $q$ is square-free and $k\in\Z$. Let $p_1,\LD, p_r$
be prime numbers $\ge 3$ such that $p_j\NDIV k$ and $\LS{q}{p_j}=1,$ $j=1,\LD,r$. Put $n=p_1p_2\cdots p_r$.
Then there are $2^r$ distinct numbers $m$, $0\le m<n$, $(m,n)=1$, such that
\rm{(\ref{2.2})} holds.

\end{corollary} 

\medskip
\noindent
{\em Proof.} Let $t$ be as in the corollary and $j\in\{1,\LD,r\}$. The congruence
\BD
  m^2-tm-1 \equiv 0 \MB{ mod } p_j
\ED
has two distinct solutions $m_1,m_2$ in $\{1,\LD,p_j-1\}$, given by
\BD
 m_1,m_2\equiv (t\pm\sqrt qk)2^* \MB{ mod } p_j,
\ED
where $\sqrt q$ denotes an integer $l$ with $\CMOD{l^2}q{p_j}$ and $2^*$ a multiplicative inverse of $2$ mod $p_j$.
Now the Chinese remainder theorem shows that the congruence
\BE
\label{2.4}
  m^2-tm-1\equiv 0 \MB{ mod } n
\EE
has $2^r$ distinct solutions $m\in\{0,\LD,n-1\}$ with $(m,n)=1$. If $m^*$ is a multiplicative inverse of $m$ mod $n$,
(\ref{2.4}) can be written $\CMOD{m-t-m^*}0n$, i.e., $\CMOD t{m-m^*}n$. Accordingly, Theorem \ref{t2} applies
to each solution $m$ of (\ref{2.4}) and gives (\ref{2.2}).
\STOP

\medskip
\noindent
{\em Remarks.} 1. If $m_1\not\equiv m_2$ mod $n$, then $1+m_1t\not\equiv 1+m_2t$ mod $nt$. So the corollary
supplies $2^r$ numbers $1+mt$ which are distinct mod $nt$ such that the corresponding \DED sums
$S(1+mt,nt)$ all have the same value. In particular, there is a set $M$ of $2^{r-1}$ numbers $1+mt$
of this kind such that all numbers in $M$ are distinct mod $nt$ and no number in $M$ has
its multiplicative inverse mod $nt$ in $M$.

2. In order to obtain examples of non-obvious equality for square-free numbers $nt$, one has to choose $t$
square-free and the primes $p_j$ such that $p_j\NDIV t$, $j=1,\LD,r$. Suppose that, in this situation,
$t$ is fixed, whereas
$r$ becomes large. Then the number $2^r$ of distinct numbers $m$ such that (\ref{2.2})
holds has the same order of magnitude as the largest possible number of arguments $m'$ for which $S(m',nt)$
can take the same value (which is $2^{r+r'}$, where $r'$ is the number of prime factors of $t$, as we
pointed out in the Introduction).

\medskip
\noindent
{\em Example.} We choose $t=7$, so $t^2+4=53$, i.e., $q=53$ and $k=1$. Further $p_1=11$, $p_2=13$, $p_3=17$
and $p_4=29$ do not divide $t$ and satisfy $\LS{53}{p_j}=1$, $j=1,\LD,4$. Hence we have $n=70499$ and
$nt=493493$. Here $m=706$ is one of 16 solutions of the congruence (\ref{2.4}). Therefore, we obtain
$16$ numbers $1+mt\in\{1,\LD,nt\}$, $(m,n)=1$, such that $S(1+mt,nt)= 2/(nt)+t/n-3\approx-2.9998966551$.
The first five of these numbers $1+mt$ are $4943$, $58535$, $79556$, $94669$, $148261$, their inverses mod $nt$ being
$488601$, $435009$, $413988$, $398875$, $345283$, respectively.

{\em Remark.} Once $n$ and $t$ have been chosen, it is possible to vary $t$. Indeed, put $t_1=t+ln$, $l\in\Z$, $l\ge 1$.
Then $\CMOD{t_1}{m-m^*}n$, so Theorem \ref{t2} also holds for $t_1$ instead of $t$. In our example, we choose
$t_1=t+2n=7+2\cdot 70499=141005$, which is the product of the primes $5$ and $28201$. Thereby, we obtain $16$
numbers $1+mt_1$ such that $S(1+mt_1,nt_1)=2/(nt_1)+t_1/n-3\approx -0.9999007076$.

\medskip
\noindent
In Corollary \ref{c4}, the crucial condition for the choice of the primes $p_j$ was
\BE
\label{2.6}
  \LS qp=1.
\EE
Whenever a prime $p\ge 3$, $p\NDIV k$,  satisfies this condition, it is eligible as one of the primes $p_j$,
$j=1,\LD,r$. It is not difficult to see that the set of primes $p\ge 3$ satisfying (\ref{2.6}) has the analytic density $1/2$
(where the set of all primes has density 1). Hence there are plenty of primes that can be chosen.
Nevertheless, it may be helpful to collect some of these primes (for small square-free  numbers $t$)
in a table (see Table 1). Since $nt$ should be square-free, we have omitted primes $p$ which divide $t$.

\BD
\begin{array}{r|r|r|l}
   t      &       q      &     k    & \hspace{9mm}      p            \\ \hline \rule[5mm]{0mm}{1mm}
   1      &       5      &     1    & 11,19,29,31,41,59 \\ \rule[4mm]{0mm}{1mm}
   2      &       2      &     2    & 7,17,23,31,41,47  \\ \rule[4mm]{0mm}{1mm}
   3      &       13     &     1    & 17,23,29,43,53,61  \\ \rule[4mm]{0mm}{1mm}
   5      &       29     &     1    & 7,13,23,53,59, 67  \\ \rule[4mm]{0mm}{1mm}
   6      &       10     &     2    & 13,31,37,41,43,53  \\ \rule[4mm]{0mm}{1mm}
   7      &       53     &     1    & 11,13,17,29,37,43  \\ \rule[4mm]{0mm}{1mm}
   10     &       26     &     2    & 11,17,19,23,37,59
\end{array}
\ED
\centerline{\bf Table 1.}

\medskip
\noindent
{\em Remarks.}
1. Further examples of equal \DED sums in the square-free case can be obtained from two theorems of Rademacher (see \cite{Ra}).
Whereas one of the nontrivial cases of his Satz 13 coincides with our case $t=1$, the other nontrivial
case concerns solutions of the congruence
\BD
  m^2-m+1\equiv 0 \MB{ mod } n.
\ED
Hence it involves square-roots of $-3$ mod $n$. His Satz 14, on the other hand, involves square-roots of $3$ mod $n$.
Therefore, this case is also not covered by Corollary \ref{c4}; indeed, it is not difficult
to see that our parameter $q$ cannot be equal to $3$.

2. For square-free positive integers $n$ with at least three prime factors, non-obvious equality seems to be a
fairly common phenomenon. For instance, let $n=7\cdot 11\cdot 13\cdot 17=17017$ and $m_1$ run through $2,3,4,5,6$.
In all of these cases there are 16 values $m_2$ such that (\ref{0.4}) holds, except $m_1=6$, where we
have only 8 values $m_2$ of this kind. Moreover, non-obvious equality occurs for all of these numbers $m_1$. For $m_1=4$,
say, there are 8 numbers $m_2$ for which $S(m_2,n)$ takes the same value; for $m_1=5$ there are 10 such numbers.
But with the exception of Corollary \ref{c4} and Rademacher's results we do not know anything for certain.

\section*{3. Proof of Theorem \ref{t1}}

Let $n$ and $d$ be positive integers and $m\in \Z$, $(m,n)=1$.
Further, let $c\in\Z$, $(c,d)=1$. Suppose that $q=md-nc$ is positive.
The three-term relation of Rademacher and Dieter
connects the \DED sums $S(m,n)$ and $S(c,d)$ in the following way:
\BE
\label{3.2}
 S(m,n)=S(c,d)+S(r,q)+\frac{n}{dq}+\frac{d}{nq}+\frac{q}{nd}-3
\EE
(see, for instance, \cite[Lemma 1]{Gi}). Here $r$ is defined as follows:
Let $j, k$ be integers such that
\BE
 \label{3.4}
  -cj+dk=1.
\EE
Then
\BE
\label{3.6}
 r=-nk+mj.
\EE
We put $d=n$ and $c=m-ln$ with $l\in\Z$, $l>0$.
Hence $q=mn- n(m-ln)=ln^2>0$.
In accordance with (\ref{3.4}), we need integers $j,k$ such that
\BD
\label{3.8}
 -mj+n(lj+k)=1.
\ED
Therefore, we may choose $j=-m^*$, where $m^*$ satisfies $\CMOD{mm^*}1n$,
and $k=(1-mm^*+nlm^*)/n=(1-mm^*)/n+lm^*$. By (\ref{3.6}),
\BD
 r=-1-lnm^*.
\ED
Since $d=n$ and $\CMOD cmn$, $S(m,n)=S(c,d)$. Accordingly, (\ref{3.2}) reads
\BE
\label{3.9}
  0=S(-1-lnm^*,ln^2)+\frac{2}{ln^2}+l-3.
\EE
If we observe $S(-1-lnm^*,ln^2)=-S(1+lnm^*,ln^2)$, we have
\BE
\label{3.10}
  S(1+lnm^*,ln^2)=\frac{2}{ln^2}+l-3.
\EE
Further, we observe that the right hand side of (\ref{3.10}) does not depend on $m$, but only on $l$ and $n$.
Hence we may replace $m^*$ by $m$, which gives
\BE
\label{3.12}
 S(1+lnm,ln^2)=\frac{2}{ln^2}+l-3.
\EE
Since $-1+lnm=-(1+ln(-m))$ and $S(1+ln(-m),ln^2)=S(1+lnm,ln^2)$, we obtain
\BE
\label{3.14}
S(-1+lnm,ln^2)=-\left(\frac{2}{ln^2}+l-3\right).
\EE
If we write $d$ instead of $l$, the identities (\ref{3.12}) and (\ref{3.14}) are just what (\ref{1.4}) says.
\STOP

\medskip
\noindent
{\em Remark.} The following modifications in the proof of Theorem \ref{t1} yield a proof of Theorem \ref{t2}:
Suppose that the positive integer $t$ is such that $\CMOD t{m-m^*}n$. Hence
$t=m-m^*+ln,$ $l\in\Z$. As in the proof of Theorem \ref{t1}, we put $d=n$, but $c=m^*-ln$ and $j=-m$.
Again $S(m,n)=S(c,d)$, and instead of (\ref{3.9})
we have
\BD
0=S(-1-mt,nt)+\frac 2{nt}+\frac tn -3,
\ED
from which (\ref{2.2}) follows.



\end{document}